# Data envelopment analysis: a flexible device to assessing SMEs' performance

Amar Oukil

*Department of Operations Management and Business Statistics, Sultan Qaboos University,*

*P.O.Box 50, Al-Khod, Sultanate of Oman*

aoukil@squ.edu.om

**Abstract:** In a context of global economy, addressing SMEs' performance within a local framework appears rather a naïve approach. The key drawback of such an approach stems from its restriction to socio-economic factors that might lead to biased decisions regarding potential venues for performance improvement. In practice, the key objective of performance analysis consists in identifying benchmarks for best managerial practices with respect to resource allocation as well as production level setting. Conducting the analysis within a specific country, let it be a developing country, may be misleading. Although, the best of the class (benchmark) can be a valid reference for its peers within the same class, its status might not be preserved if the analysis is projected outside the borders of the class. Indeed, the likelihood for outperformance is high. In order to set targets for global competition, decision makers ought to look at the concept of performance from a broader geographical perspective, instead of confining it to a local scope. Here, we analyze, through a case study, SMEs' performance within local and global production technology frameworks and we highlight the impact of the economy scope on various decisions. Data envelopment analysis (DEA) is used as a mathematical tool to support such decisions.

**Introduction**

Small- and medium sized enterprises (SME's) are important economic entities (Tagliavini et al. 2001), often associated with nations' growth (Beck et al. 2005) and having a significant impact on employment (Reynolds 1997). Therefore, serving effectively this category of firms appears a real priority for public business services (Shi et al. 2013). Primarily, such a task requires a reliable evaluation of the potential of these enterprises for prospective improvement through efficient utilization of on-hand resources (Satayapaisal et al. 2012). Performance of resource allocation inside the enterprise can be assessed over various economic scopes using data envelopment analysis (DEA) (Cooper et al. 2002). DEA is a non-parametric approach for the estimation of efficiency ratios as well as the input reduction or/and output expansion required from inefficient units in order to reach full efficiency (Soltani *et al*. 2021). In practice, DEA goes far beyond the computation of the efficiency; typically, it allows the decision maker to know what operating practices, mix of resources, scale sizes, scope of activities, and more, the operating units may adopt to improve their performance (Oukil and Govindaluri, 2020). Several studies pertaining to the application of DEA to assessing SMEs performance can be found in the literature. For instance, Rani et al. (2014) used DEA within a simulation framework to identify the problems faced by SMEs in the Malaysian food industry and set a template of actions. Průša (2012) applied DEA to perform structural analysis of Czech SMEs in manufacturing based on their relative efficiency. Roveicy et al. (2015) also employed DEA to compare the performance of knowledge based firms (R&D SMEs) across different provinces of Iran. Monelos et al. (2014) employed DEA as a proxy approach, jointly with discriminant, logit and multivariate linear models, to diagnose and forecast business failure for non-financial SMEs in Galicia. In the tofu industry, Purwanto et al. (2014) adopted DEA models to estimate the ability of this sector's SMEs to survive and compete in Indonesia. Apparently, most of the proposed studies are restricted to specific countries. Although the findings of related investigations are acceptable country wise, there is doubt for such a status to be preserved with reference to a more extended worldwide set. The key drawback of such an approach stems from its restriction to socio-economic factors that might lead to biased decisions. Conducting the analysis within a specific country, let it be a developing country, may be misleading with respect to most managerial aspects, including identifying benchmarks for best practice, resource reallocation as well as production level setting. Although, the best of the class (benchmark) can be a valid reference for its peers within the same class, its status might be downgraded should the analysis be projected outside the borders of the class. Indeed, the relative aspect of performance analyses raises in-depth issues regarding decision consistency in a rolling space. In order to set targets for global competition, decision makers ought to look at the concept of performance from a broader geographical perspective, instead of confining it to a local scope. In spite of its vital importance for many businesses, this problem has never been addressed explicitly in the literature. Here, we analyze, through a case study, SMEs' performance within local and global production technology frameworks and we highlight the impact of the economy scope on various key decisions. Our paper unfolds as follows. We first explain the methodological context of the current study, focusing on the key DEA concepts. Next, we present succinctly the contextual setting alongside the results of the application of the DEA model. We conclude with a discussion and recommendations for future research.

**Methodological background**

Data envelopment analysis (DEA) is a non-parametric approach for the estimation of relative efficiency ratios of decision making units (DMUs). The efficiencies assessed reflect the scope for resource conservation at the unit being assessed without detriment to its outputs or, alternatively, the scope for output augmentation without additional resources (Al-Mezeini *et al.* 2021). Under the assumption of constant returns to scale (CRS), an increase in all inputs at a certain rate will cause an increase in output at the same rate. However, if the DMUs are assessed within the variable returns to scale (VRS) context, the change between inputs and outputs will not occur at the same rate. Figure 1 shows an example of six DMUs, A to F, evaluated under CRS assumption. The production possibility set (PPS), a.k.a. the efficiency frontier, is the line OA, which depicts A as the only efficient DMU. All the remaining DMUs are deemed to be inefficient.

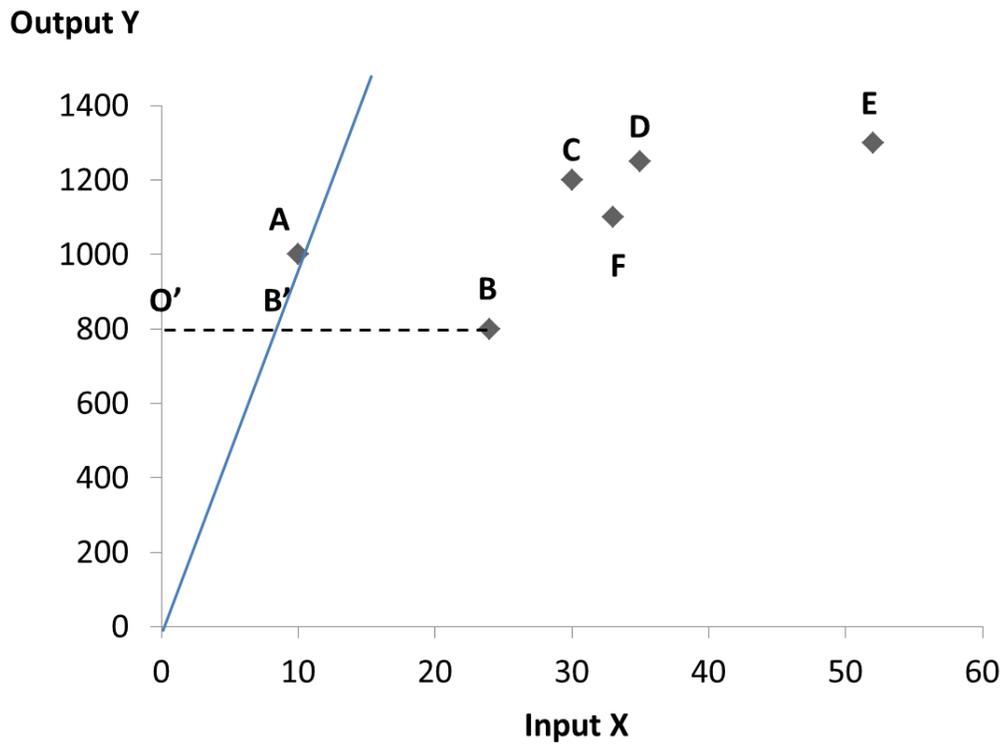

**Figure 1**: Production possibility set under CRS assumption

The same set of DMUs is represented in Figure 2 under VRS assumption. Here, the PPS is a piecewise curve (Sow *et al*. 2016), where only two DMUs, namely B and F, are inefficient.

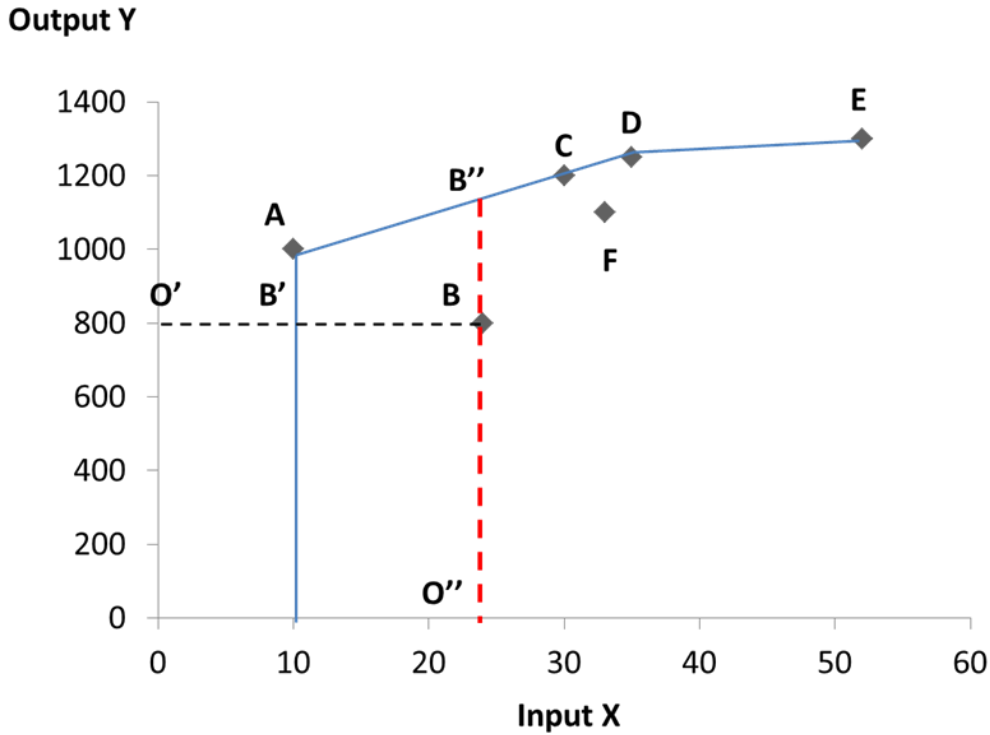

**Figure 2**: Production possibility set under VRS assumption

The computation of the relative efficiency θ of a given DMU may be based on either input reduction or output expansion. If an input oriented DEA model is used, θ is the proportion to which the observed input levels can be reduced for given output levels (Amin and Oukil, 2019). On the other hand, if the model is output oriented, θ becomes the proportion of observed output levels to maximum possible output levels for given input levels. In Figure 1, A being part of the PPS, its relative efficiency is 100%. The efficiency of B, for instance, is computed as the ratio of target input O'B' over the actual input O'B on the input axis, i.e., $\theta_B$= O'B'/O'B. On the output axis, even if $\theta_B$ remains the same under the CRS assumption, its value increases under VRS (Oukil and El-Bouri, 2021).

To guarantee a better focus on the objectives initially set to our study, all the investigations will be conducted through input oriented models under CRS assumption. In what follows, the DMU is an SME.

Assume a set of $K$ SMEs, with each SMU $k$ being defined with $N$ inputs $x$ and $M$ outputs $y$. In reference to the underlying production technology, SME $(x_k, y_k)$ is fully defined with the observed values of $x_{ik}$ and $y_{jk}$, $i=1,.., N$ and $j=1,.., M$. To estimate the efficiency score $\theta_f$ of SME $(x_f, y_f)$ and set production targets for inefficient SMEs, the input-oriented formulation of the CCR model (Charnes, Cooper, and Rhodes 1978), can be written as follows.

$$\min \theta_f \qquad (1)$$

Subject to:

(CCR)
$$\sum_{k=1}^{K} \lambda_k x_{ik} \leq \theta_h x_{ih} \quad i = 1,...,N \qquad (2)$$

$$\sum_{k=1}^{K} \lambda_k y_{jk} \geq y_{jh} \quad j = 1,...,M \qquad (3)$$

$$\lambda_k \geq 0 \quad k = 1,...,K \qquad (4)$$

The efficiency $\theta_f$ of SME $(x_f, y_f)$ represents the minimum radial decrease of inputs that is required to reach the efficiency frontier for a specified level of inputs. The vector $\lambda$ measures the weights of peers in producing the projection of SME $(x_f, y_f)$ on the efficiency frontier (Oukil and Al-Zidi, 2018). Constraints (2) and (3) state that reference points are linear combinations of the input and output values of efficient peers for SME $(x_f, y_f)$

BCC model (Banker, Charnes, and Cooper 1984) can be obtained from (CCR) by adding the convexity constraint guaranteeing that only weighted averages of efficient SMEs enter the reference set, i.e.,

$$\sum_{k=1}^{K} \lambda_k = 1 \qquad (5)$$

CCR and BCC models are both formulated with the implicit assumption that the assessed SMEs operate within homogeneous environments, which presupposes that only variables representing proper inputs are an integral part of the production technology (Hassan and Oukil, 2021). The optimal solutions of CCR and BCC models, denoted $\theta^*_{CCR}$ and $\theta^*_{BCC}$ respectively, are known as aggregate and technical efficiency scores of SME $(x_f, y_f)$.

**Application and results**

A sample of 1000 SMEs split among 9 countries, A to I, has been generated through 10 input and 6 output variables, as shown in Table 1. The number of inputs and outputs has been selected in a way that guarantees clear efficiency discrimination within each country, i.e., the number of SMEs $\Psi$ satisfies the inequality $\Psi \geq \max[MN, 3 \times (M+N)]$ (Cooper *et al.* 2002).

**Table 1**: Distribution of the SMEs over the nine countries

| Country | A | B | C | D | E | F | G | H | I |
|---|---|---|---|---|---|---|---|---|---|
| **Number of SMEs** | 68 | 136 | 90 | 100 | 120 | 140 | 70 | 111 | 165 |

The detailed sets of data are available from the author under request.

We use IBM-ILOG CPLEX version 12.4 to compute the optimal efficiency scores $\theta^*$ and the slack values for each SME within each country, then on a global scale, where all the SMEs are grouped in the same cluster.

Table 2: Country-wise efficiency results of the SMEs over the nine countries

| Country | A | B | C | D | E | F | G | H | I |
|---|---|---|---|---|---|---|---|---|---|
| Avg | 0.972 | 0.971 | 0.974 | 0.975 | 0.986 | 0.963 | 0.987 | 0.976 | 0.983 |
| STD | 0.078 | 0.083 | 0.102 | 0.084 | 0.050 | 0.106 | 0.050 | 0.084 | 0.064 |
| min | 0.677 | 0.572 | 0.365 | 0.590 | 0.701 | 0.255 | 0.657 | 0.362 | 0.570 |
| Efficient | 58 | 115 | 82 | 84 | 104 | 115 | 63 | 98 | 145 |
| % | 85.29 | 84.56 | 91.11 | 84.00 | 86.67 | 82.14 | 90.00 | 88.29 | 87.88 |

Table 2 provides the results of the efficiency evaluation of the SMEs within each country. The first row shows relatively a high level of efficiency, with averages 96% for all the countries. Furthermore, the proportions of efficient SMEs, which are also the potential benchmarks for inefficient SMEs, ranges between 82.14% and 91.11%. Although such figures may suggest that the SMEs are, on average, performing well within their respective countries, the assessment remains locally bounded.

In order to detect the source of inefficiency for non-efficient SMEs, we conduct a slack analysis on the ten inputs for each country and we calculate the corresponding average excess input proportion. The results are given in Table 3.

Table 3: Country-wise slack analysis of the SMEs over the nine countries

| Country | A | B | C | D | E | F | G | H | I |
|---|---|---|---|---|---|---|---|---|---|
| X1 | 0.54 | 1.07 | 2.18 | 3.26 | 2.27 | 1.95 | 3.12 | 2.55 | 1.02 |
| X2 | 1.54 | 2.11 | 2.00 | 2.50 | 1.88 | 2.11 | 2.31 | 2.44 | 2.52 |
| X3 | 2.12 | **2.92** | 0.61 | 1.84 | 2.86 | 1.12 | 0.29 | 1.20 | 2.80 |
| X4 | 2.04 | 1.63 | 1.95 | 5.12 | 1.78 | 2.87 | 4.05 | 1.79 | 1.12 |
| X5 | 2.96 | 2.58 | 1.15 | 5.05 | 2.53 | 0.72 | 3.19 | 1.05 | 1.83 |
| X6 | 3.25 | 4.04 | 1.18 | 1.22 | 2.22 | **4.87** | **4.75** | **2.69** | **2.81** |
| X7 | **6.68** | 2.56 | 2.04 | **5.95** | 1.77 | 3.66 | 3.56 | 2.64 | 1.82 |
| X8 | 2.44 | 0.46 | 1.95 | 5.69 | **3.62** | 1.18 | 0.12 | 1.26 | 1.95 |
| X9 | 1.80 | 2.28 | **2.10** | 4.15 | 2.58 | 2.41 | 1.64 | 1.29 | 1.41 |
| X10 | 1.16 | 2.83 | 1.65 | 3.85 | 1.72 | 3.34 | 1.20 | 1.43 | 2.29 |

The largest proportions of excess input used throughout the nine countries range between 2.10% and 6.68%, regardless of the resource that the input refers to. This range indicates that most of the inefficient SMEs are able to improve their performance without need to too much resource reduction, that is, most of these SMEs are very close to their local PPS.

In a second stage, similar computations are performed for a sample grouping all 1000 SMEs. The disaggregated values obtained for each country using the large sample's results are summarized in Tables 4 and 5.

Table 4: Global-wise efficiency results of the SMEs over the nine countries

| Country | A | B | C | D | E | F | G | H | I |
|---|---|---|---|---|---|---|---|---|---|
| Avg | 0.876 | 0.887 | 0.915 | 0.878 | 0.923 | 0.900 | 0.911 | 0.876 | 0.894 |
| STD | 0.182 | 0.167 | 0.173 | 0.168 | 0.133 | 0.161 | 0.151 | 0.167 | 0.151 |
| min | 0.387 | 0.395 | 0.263 | 0.311 | 0.481 | 0.241 | 0.420 | 0.277 | 0.458 |
| Efficient | 37 | 75 | 66 | 53 | 76 | 85 | 45 | 61 | 88 |
| % | 54.41 | 55.15 | 73.33 | 53.00 | 63.33 | 60.71 | 64.29 | 54.95 | 53.33 |
| Worse | 31 | 61 | 24 | 47 | 44 | 55 | 25 | 50 | 77 |
| % | 45.59 | 44.85 | 26.67 | 47.00 | 36.67 | 39.29 | 35.71 | 45.05 | 46.67 |
| Shifted | 21 | 41 | 16 | 31 | 28 | 30 | 18 | 37 | 60 |
| % | 30.88 | 30.15 | 17.78 | 31.00 | 23.33 | 21.43 | 25.71 | 33.33 | 36.36 |

Global-wise, there is clearly a decline in the number of efficient SMEs within all the countries, with the proportions now ranging from 53% and 73.33%. Overall, 41.4% of the SMEs have their performance deteriorated, as compared to the local assessment (country-wise). In the meantime, the proportions of SMEs whose efficiency scores have worsen when considered globally fall between 26.67% and 47%, as shown in Table 4. Among the latter worsening SMEs, 17.78% to 36.36% have shifted status from efficient to inefficient, that is, they lost their benchmarking advantage.

Table 5: Global-wise slack analysis of the SMEs over the nine countries

| Country | A | B | C | D | E | F | G | H | I |
|---|---|---|---|---|---|---|---|---|---|
| X1 | 0.96 | 2.92 | 2.73 | 4.88 | 3.40 | 3.51 | 6.31 | 3.65 | 5.02 |
| X2 | 3.32 | 4.68 | 2.51 | 4.44 | 4.87 | 3.76 | 3.58 | 4.62 | 4.44 |
| X3 | 3.47 | **7.45** | 1.78 | 4.60 | 4.64 | 4.62 | 4.78 | 2.80 | 5.27 |
| X4 | **5.06** | 2.44 | 2.61 | 3.92 | 3.98 | 3.06 | 4.94 | 4.41 | 4.12 |
| X5 | 4.43 | 5.00 | 1.83 | **6.28** | 4.88 | 2.27 | 3.02 | **5.53** | **6.03** |
| X6 | 3.27 | 4.21 | 2.06 | 1.79 | 5.92 | 1.90 | **7.84** | 4.64 | 3.60 |
| X7 | 3.52 | 2.83 | **3.11** | 2.04 | 4.78 | 4.00 | 3.06 | 2.65 | 3.78 |
| X8 | 0.70 | 3.80 | 2.14 | 2.88 | 2.43 | 2.66 | 2.07 | 1.72 | 4.43 |
| X9 | 2.36 | 3.74 | 2.74 | 4.92 | 5.89 | **4.65** | 4.53 | 3.74 | 3.72 |
| X10 | 4.49 | 5.18 | 1.10 | 2.90 | **7.54** | 3.97 | 3.92 | 4.75 | 5.26 |
| Status | X | - | X | X | X | X | - | X | X |

The effect of the global approach over the performance evaluation process can also be gauged via the slack results shown in Table 5. It is obvious that the excess inputs used have increased for most resources as a direct result of the decline in the performance of the majority of SMEs. On another hand, the pivotal input for performance improvement has changed for most of the countries, except for countries B and G.

**Conclusion and recommendations**

"In the land of the blind, the one-eyed man is king"
Although the above quote has been stated for individuals, it appears strongly well-fitted to the context of our study. We departed with a performance analysis, performed locally, where more than 82% of the SMEs are efficient over the nine countries. Practically speaking, more than 82% of the SMEs are evaluated as the best of the class and, as a consequence, they can be seen as the benchmarks for the least performing SMEs. In a closed universe, like a country, this outcome may raise some SMEs to the leading positions, kings in their world, without being really as excellent as evaluated.

In a global economy, ambitious SMEs ought to look at their performance beyond the local geographical borders. Such a move will not only allow the SME to measure its potential according to real benchmarks but and foremost to propel the quality of its management towards international standards.

The work presented in this paper is no more than an introduction to a topic that deserves to be handled with real-life data in a view to develop indicators to measure the consistency of various decisions. Data envelopment analysis, with its mathematical properties, is undoubtedly one of the most appropriate tools to conducting such investigations successfully. Further development may consider the impact of the contextual setting on the ranking of the SMEs, To this end, one could consider more advanced DEA tools, like DEA cross-efficiency (e.g., Oukil, 2020).